\documentclass[12pt,a4paper]{article}

\usepackage{graphics}
\usepackage[dvips]{graphicx}
\usepackage{t1enc}
\usepackage{amscd}
\usepackage{amsmath,amssymb}
\usepackage{theorem,url}
\usepackage{subfig}
\usepackage[version=3]{mhchem}

\renewcommand{\emph}{\textbf}
\newcommand{\mma}[1]{\textbf{\texttt{#1}}}

\theorembodyfont{\upshape}
\newtheorem{Thm}{Theorem}
\newtheorem{Def}{Definition}

\newtheorem{Lem}{Lemma}

\newcommand{\torol}[1]{}
\newcommand{\nulb}{\ifmmode \mathbf{0}\else \textbf{0}\fi}

\newcommand{\Rp}{\mathbb{R}^+}

\newcommand{\Mma}{\textit{Mathematica }}

\DeclareMathOperator{\rank}{rank}

\DeclareMathOperator{\Span}{span}
\DeclareMathOperator{\diag}{diag}
\DeclareMathOperator{\Ran}{Ran}

\title{Global Controllability of Chemical Reactions}
\author{
D\'aniel Andr\'as Drexler\\
Department of Control Engineering and Information Technology\\
Budapest University of Technology and Economics,\\
Budapest, Magyar tud\'osok krt. 2., HUNGARY, H-1117\\
\url{drexler@iit.bme.hu} and\\
J\'anos T\'oth\\
Department of Mathematical Analysis\\
Budapest University of Technology and Economics,\\
Budapest, Egry J. u. 1., HUNGARY, H-1111 and\\
Laboratory for Chemical Kinetics
Eötvös University\\
Budapest, Pázmány Péter sétány 1/A, HUNGARY, H-1117\\
\url{jtoth@math.bme.hu}}
\date{}

\begin{document}
{\LARGE
\noindent \textbf{Title Page}

\noindent \textbf{Title:}
Global Controllability of Chemical Reactions

\noindent \textbf{Authors:}

D\'aniel Andr\'as Drexler\footnote{(author designated to review proofs)
Department of Control Engineering and Information Technology\\
Budapest University of Technology and Economics,\\
Budapest, Magyar tud\'osok krt. 2., HUNGARY, H-1117\\
\url{drexler@iit.bme.hu}} and\\

J\'anos T\'oth\footnote{
Department of Mathematical Analysis\\
Budapest University of Technology and Economics,\\
Budapest, Egry J. u. 1., HUNGARY, H-1111 and\\
Laboratory for Chemical Kinetics
Eötvös Loránd University\\
Budapest, Pázmány Péter sétány 1/A, HUNGARY, H-1117\\
\url{jtoth@math.bme.hu}}

\noindent \textbf{Running Head}: Global Controllability of Chemical Reactions}
\vfill\eject

\maketitle

\begin{abstract}
Controllability of chemical reactions is an important problem in chemical engineering science. In control theory, analysis of the controllability of linear systems is well-founded, however the dynamics of chemical reactions is usually nonlinear. Global controllability properties of chemical reactions are analyzed here based on the Lie-algebra of the vector fields associated to elementary reactions. A chemical reaction is controllable almost everywhere if all the reaction rate coefficients can be used as control inputs. The problem where one can not control all the reaction rate coefficients is also analyzed. The reaction steps whose reaction rate coefficients need to be control inputs are identified. A general definition of consecutive reactions is given, and it turns out that they are controllable almost everywhere by using the reaction rate coefficient of only one reaction step as control input.
\end{abstract}

{\bf Keywords}: kinetics, chemical reactions, controllability, Lie-algebra
\section{Introduction}

Controlling chemical reactions is a key issue in modern chemical engineering science, see e.g. \cite{chemicalreactors}, where control of nonlinear chemical reactors is considered. The well founded methods in control engineering are based on linear dynamical systems, however the dynamics of chemical reactions is usually nonlinear. Control and controllability of nonlinear systems is only available for a relatively small class of systems, and is more involved than control of systems with linear dynamics, see e.g. \cite{isidori,coron}. Controllability is a fundamental property of dynamical systems, since it tells us whether it is possible to reach the desired goal or not.

Controllability of chemical reactions is usually analyzed using control theory developed for linear systems, and the linear model is acquired by linearizing the dynamics at an operating point. Controllability of chemical reactions that have a positive equilibrium is established in \cite{farkascontrol} based on the linearized dynamics. In \cite{catalyticoxidation}, controllability analysis of a liquid-phase catalytic oxidation of toluene to benzoic acid was done based on linearization at five different operating points, and controllability analysis of polymerization in different operating points was done in \cite{polymerization}. Controllability analysis of protein glycosylation was done based on a linear model identified from measurements in \cite{glycosylation}. The connection of controllability and structure of chemical reaction was analyzed in \cite{szederkenyi2008}, while controllability of chemical processes was analyzed in \cite{chemicalprocess,chemicalprocess2,designandcontrol,szederkenyi2002nonlinear,otero2006dynamic}. 

Controllability analysis in operating points however only gives local results. The application of Lie-algebra rank condition was suggested for the controllability analysis of chemical reactions in \cite{ervadiradhnakrishnanvoit}, however no general results were given. A review of controllability analysis of chemical reactions was done in \cite{overview2011}.

We analyze the Lie-algebra of the vector fields related to the reaction steps, and use the Lie-algebra rank condition as suggested in \cite{ervadiradhnakrishnanvoit} to get global controllability results for chemical reactions. We consider the reaction rate coefficients of the reaction steps as control inputs and seek to find the lowest number of control inputs needed to control the system. The analysis done in this paper is symbolical, contrary to the analysis done in the literature that is usually numerical with some few exceptions, e.g. \cite{farkascontrol}. Our investigation will be restricted to those reactions for which the stoichiometric and kinetic subspaces coincide \cite{FeinbergHorn}; however this is not a too strong restriction.

We give the model structure of kinetic equations that is used to model chemical reactions, and the Lie-algebra rank condition in Section \ref{sec:preliminaries}. We show that chemical reactions are controllable almost everywhere if the control inputs are the reaction rate coefficients of the reaction steps in Section \ref{sec:allinput}. In Section \ref{sec:lessinput}, we further analyze the Lie-algebra of the vector fields corresponding to reaction steps, and give results on controllability when not all the reaction rate coefficients are control inputs. We give a lower and upper bound for the required number of control inputs, and identify the reaction steps whose reaction rate coefficients need to be control inputs. We define the initializer reaction steps and prove that their reaction rate coefficients need to be control inputs, that has already been shown by experiments in \cite{polymerization} on a polymerization example. We define the consecutive reactions, and show that they can be controlled with only one control input. We give two examples for consecutive reactions in Section \ref{sec:examples}, and analyze their controllability and region of controllability. The paper ends with a discussion in Section \ref{sec:discuss}.

\section{Kinetics of reactions and controllability}\label{sec:preliminaries}

\subsection{The differential equation of reaction kinetics}
Let the number of species participating in a chemical reaction be $M$, and the species be denoted by $\ce{X1},\ce{X2},\ldots,\ce{X_M}$. Suppose, that the reaction consists of $R$ reaction steps. We call the linear combination of the species on the left-hand side of a reaction step the reactant complex, and the linear combination of the species on the right-hand side of a reaction step the product complex, so a reaction step is defined by giving the reactant complex that transforms into a product complex with a given reaction rate coefficient.  Let $\alpha(m,r)$ denote the number of species $\ce{X_m}$ in the reactant complex of the $r$th reaction step and let $\beta(m,r)$ denote the number of species $\ce{X_m}$ in the product complex of the $r$th reaction step. Let $k_r$ be the reaction rate coefficient of the $r$th reaction step. The chemical reaction system (or simply: reaction) is defined by the $R$ reaction steps \cite{erditoth}
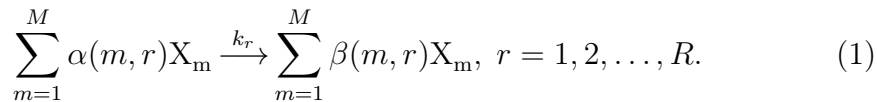
\begin{equation}\label{eq:kemreak1}
\mathop{\sum}_{m=1}^M \alpha(m,r)\ce{X_m} \mathop{\longrightarrow}^{k_r} \mathop{\sum}_{m=1}^M \beta(m,r)\ce{X_m}, \; r=1,2,\ldots,R.
\end{equation}
The matrix $\boldsymbol{\gamma} = \boldsymbol{\beta}-\boldsymbol{\alpha}$ is called the stoichiometric matrix, and its columns are called the reaction step vectors. The usual mass action type differential equation of the reaction (\ref{eq:kemreak1}) is
\begin{equation}\label{eq:kemdiff1}
\dot{\boldsymbol{x}}(t) = \mathop{\sum}_{r=1}^R k_r\boldsymbol{\gamma}(\cdot,r)\boldsymbol{x}(t)^{\boldsymbol{\alpha}(\cdot,r)}
\end{equation}
where $x_m(t)$ is the concentration of the species $\ce{X_m}$ at time $t$, $\boldsymbol{x}(t)$ is the vector of concentrations with its $m$th element being $x_m(t)$, $\boldsymbol{\gamma}(\cdot,r)$ denotes the $r$th column of the matrix $\boldsymbol{\gamma}$, and $\boldsymbol{x}^{\alpha(\cdot,r)}$ is the monomial
\begin{equation}
\boldsymbol{x}^{\alpha(\cdot,r)} = x_1^{\alpha(1,r)}x_2^{\alpha(2,r)} \cdot \ldots \cdot x_M^{\alpha(M,r)}.
\end{equation}
\begin{Def}
We call the image space of the stoichiometric matrix $\boldsymbol{\gamma}$ the stoichiometric space.
\end{Def}
The stoichiometric space is the tangent space of the system (\ref{eq:kemdiff1}) in almost all practically relevant cases. As to the exceptions, which we disregard here, see \cite{FeinbergHorn}. Since the stoichiometric space is the tangent space of the system (\ref{eq:kemdiff1}), the dynamical system is only able to change its state locally in the sets of the form $(\boldsymbol{x}_0 + \Ran \boldsymbol{\gamma})\cap (\Rp_0)^M$, called reaction simplexes \cite{erditoth}.

\subsection{Controllability of nonlinear systems and the Lie-algebra rank condition}

Consider an input affine system with the differential equation and initial condition
\begin{equation}\label{eq:nonlin1}
\dot{\boldsymbol{x}}(t)=\boldsymbol{f}(\boldsymbol{x}(t)) + \mathop{\sum}_{j=1}^J \boldsymbol{g}_j(\boldsymbol{x}(t))u_j(t),\quad \boldsymbol{x}(0)=\boldsymbol{x}_0 \in \mathbb{R}^M.
\end{equation}
We say that the system is smooth if $\boldsymbol{f}$ and $\boldsymbol{g}_1,\boldsymbol{g}_2,\ldots,\boldsymbol{g}_J$ are smooth functions, i.e. $\boldsymbol{f},\boldsymbol{g}_1,\boldsymbol{g}_2,\ldots,\boldsymbol{g}_J \in \mathcal{C}^\infty (\mathbb{R}^M,\mathbb{R}^M)$.

\begin{Def}
System (\ref{eq:nonlin1}) is called controllable at a point $\boldsymbol{x}^* \in \mathbb{R}^M$ if there exists an open neighborhood $U$ of $\boldsymbol{x}^*$ such that if $\boldsymbol{x}(0) \in U$ then there exists $u_1(t), u_2(t), \ldots , u_J(t),  t\geq 0$ such that $\mathop{\lim}_{t\rightarrow +\infty} \boldsymbol{x}(t) = \boldsymbol{x}^*$.  
\end{Def}

Note that we only deal with controllability, and do not discuss the problem of finding the input functions that are required to get into the desired state.

A conventional way of analyzing controllability of nonlinear smooth input affine systems like (\ref{eq:nonlin1}) is based on the Lie-algebra generated by the control vector fields $\boldsymbol{g}_j$, $j=1,2,\ldots,J$ and the Lie-bracket of the drift vector field $\boldsymbol{f}$ and the control vector fields. 

\begin{Def}
The Lie-bracket of two smooth vector fields $\boldsymbol{\varphi},\boldsymbol{\psi} \in \mathcal{C}^\infty(\mathbb{R}^M,\mathbb{R}^M)$ is
\begin{equation}\label{eq:liebracket}
[\boldsymbol{\varphi},\boldsymbol{\psi}] = \boldsymbol{\psi}'\boldsymbol{\varphi}-\boldsymbol{\varphi}'\boldsymbol{\psi}.
\end{equation}
\end{Def}

\begin{Def}
A Lie-algebra \cite{Humphreys} is a vector space $V$ over the field $\mathbb{F}$ equipped with the binary operation $[\cdot,\cdot]$ called the Lie-bracket (or Lie-product), if it satisfies the following properties.
\begin{enumerate}
\item Bilinearity, i.e. $\forall c_1,c_2 \in \mathbb{F}$ and $\forall \boldsymbol{x},\boldsymbol{y},\boldsymbol{z} \in V$ 
\begin{eqnarray}
[c_1\boldsymbol{x}+c_2\boldsymbol{y},\boldsymbol{z}]&=&c_1[\boldsymbol{x},\boldsymbol{z}]+c_2[\boldsymbol{y},\boldsymbol{z}]\\
\left[ \boldsymbol{z} , c_1 \boldsymbol{x} + c_2 \boldsymbol{y} \right] &=& c_1 [ \boldsymbol{z} , \boldsymbol{x} ]+ c_2 [ \boldsymbol{z} , \boldsymbol{y} ]
\end{eqnarray}
\item $\forall \boldsymbol{x} \in V$, $[\boldsymbol{x},\boldsymbol{x}]=0$.
\item The Jacobian identity, i.e. $\forall \boldsymbol{x},\boldsymbol{y},\boldsymbol{z} \in V$
\begin{equation}
[\boldsymbol{x},[\boldsymbol{y},\boldsymbol{z}]]+[\boldsymbol{y},[\boldsymbol{z},\boldsymbol{x}]]+[\boldsymbol{z},[\boldsymbol{x},\boldsymbol{y}]]=0.
\end{equation}
\end{enumerate}
\end{Def}
The first two properties of the Lie-bracket imply anticommutativity, i.e. $\forall \boldsymbol{x},\boldsymbol{y} \in V$, $[\boldsymbol{x},\boldsymbol{y}]=-[\boldsymbol{y},\boldsymbol{x}]$. In this article, we will use the Lie-algebra on the vector space of smooth functions $\mathcal{C}^\infty (\mathbb{R}^M,\mathbb{R}^M)$ over the field $\mathbb{R}$ equipped with the Lie-bracket defined by (\ref{eq:liebracket}).

\begin{Def}
The Lie-algebra $\mathcal{L}$ generated by the vector fields $\boldsymbol{f}_1$ and $\boldsymbol{f}_2$, denoted by $\mathcal{L}=Lie(\boldsymbol{f}_1,\boldsymbol{f}_2)$ is the smallest subspace of $\mathcal{C}^\infty(\mathbb{R}^M,\mathbb{R}^M)$ that satisfies the following properties \cite{coron}:
\begin{enumerate}
\item $\boldsymbol{f}_1,\boldsymbol{f}_2 \in \mathcal{L}$;
\item If $\boldsymbol{x},\boldsymbol{y} \in \mathcal{L}$, then $[\boldsymbol{x},\boldsymbol{y}] \in \mathcal{L}$.
\end{enumerate}
\end{Def}

\begin{Def}
The object that assigns to every point $\boldsymbol{x} \in \mathbb{R^M}$ a linear subspace of $\mathbb{R}^M$ is called a distribution \cite{isidori}. Distributions are usually denoted by the symbol $\Delta$.
\end{Def}

\begin{Def}
A vector field $\boldsymbol{f}$ is said to be an element of a distribution $\Delta$ ($\boldsymbol{f} \in \Delta$), if the value of $\boldsymbol{f}$ is in the subspace assigned to $\Delta$ in every point, i.e. $\forall \boldsymbol{x} \in \mathbb{R}^M$, $\boldsymbol{f}(x) \in \Delta(x)$.
\end{Def}

\begin{Def}
A distribution $\Delta_C$ is called the controllability distribution \cite{isidori} of the system (\ref{eq:nonlin1}), if
\begin{enumerate}
\item $\forall \boldsymbol{x} \in  \mathbb{R}^M$, $\Span \left\{ \boldsymbol{g}_1(\boldsymbol{x}),\boldsymbol{g}_2(\boldsymbol{x}),\ldots,\boldsymbol{g}_J(\boldsymbol{x}) \right\} \subseteq \Delta_C$.
\item It is invariant under the vector field $\boldsymbol{f}$, i.e. if $\boldsymbol{\tau}\in \Delta_C$, then $[\boldsymbol{\tau},\boldsymbol{f}] \in \Delta_C$.
\item It is involutive, i.e. if $\boldsymbol{\tau}_1,\boldsymbol{\tau}_2 \in \Delta_C$, then $[\boldsymbol{\tau}_1,\boldsymbol{\tau}_2] \in \Delta_C$.
\end{enumerate}
\end{Def}
The above conditions mean that the controllability distribution is spanned by the vector fields from the Lie-algebra generated by the vector fields $\boldsymbol{g}_j$ and $[\boldsymbol{f},\boldsymbol{g}_j]$, $j=1,2,\ldots,J$ \cite{coron}.

\begin{Def}
The controllable subspace of the system (\ref{eq:nonlin1}) is the linear combination of all the space variable vectors that can be reached using appropriate inputs $u_1,u_2,\ldots,u_J$.
\end{Def}

\begin{Thm}
The dimension of the subspace $\Delta_C(\boldsymbol{x}^*)$ is the dimension of the controllable subspace of the system at the point $\boldsymbol{x}^*$, so the system is controllable in $\boldsymbol{x}^*$ if and only if $\dim \Delta_C(\boldsymbol{x}^*) = M$.
\end{Thm}
This is also called the Lie-algebra rank condition \cite{isidori, coron}. Note that for linear systems with the differential equation
\begin{equation}\label{eq:linearmodel}
\dot{\boldsymbol{x}}(t)=\boldsymbol{A}\boldsymbol{x}(t)+\boldsymbol{B}u(t)
\end{equation}
the controllability distribution assigns the same subspace for every point, and this subspace is the image space of the controllability matrix
\begin{equation}\label{eq:controlmatrix}
\boldsymbol{M}_C = \left[\begin{array}{ccccc} \boldsymbol{B} & \boldsymbol{A}\boldsymbol{B} & \boldsymbol{A}^2\boldsymbol{B} & \ldots & \boldsymbol{A}^{M-1}\boldsymbol{B} \end{array}\right].
\end{equation}
So a linear system is controllable if and only if $\rank \boldsymbol{M}_C=M$. For example, in \Mma\, the function \mma{ControllableModelQ} can be used to determine if a linear system is controllable or not.

\section{Controllability of chemical reactions if all reaction rate coefficients are control inputs}\label{sec:allinput}
As a first step we consider chemical reactions in the form of (\ref{eq:kemdiff1}) assuming that all of the reaction rate coefficients can be varied independently, so all reaction rate coefficients are control inputs. The chemical reaction under consideration thus is the driftless smooth input affine system
\begin{equation}\label{eq:kemdiff2}
\dot{\boldsymbol{x}}(t) = \mathop{\sum}_{r=1}^R \boldsymbol{\gamma}(\cdot,r)\boldsymbol{x}(t)^{\boldsymbol{\alpha}(\cdot,r)}u_r(t)
\end{equation}
that is obtained from (\ref{eq:kemdiff1}) with the substitution $k_r = u_r(t)$ for $r = 1,2,\ldots,R$. This chemical reaction consists of $R$ reaction steps with $M$ species. Usually the $M$ differential equations in (\ref{eq:kemdiff2}) are not linearly independent, in spite of the fact that usually $R$ is much larger than $M$; in combustion chemistry, the Law's law \cite{Law} says that $R \approx 5 M$.

Since the trajectories of a chemical reaction system are in the reaction simplex that is contained in $\Ran \boldsymbol{\gamma}$, the dimension of the controllable subspace can not be larger than the dimension of the image space of $\boldsymbol{\gamma}$.
\begin{Def}
We say that a chemical reaction is controllable at a point $\boldsymbol{x}^* \in \mathbb{R}^M$ if $\dim \Delta_C(x^*)=\rank \boldsymbol{\gamma}$. 
\end{Def}
The controllability of chemical reactions in the form of (\ref{eq:kemdiff2}) have already been analyzed in \cite{farkascontrol}, where the author proved the following two theorems:
\begin{Thm} 
Suppose that the system (\ref{eq:kemdiff2}) has a positive equilibrium $\boldsymbol{x}^*$ for $u_1=u_2=\ldots=u_R=1$. Then for the control matrix $\boldsymbol{B}$ of the linear system (\ref{eq:linearmodel}) obtained from (\ref{eq:kemdiff2}) by linearization at the equilibrium point $\boldsymbol{x}^*$, we have that $\rank \boldsymbol{B} = \rank \boldsymbol{\gamma}$, i.e. the linearized model is controllable, meaning that the system is controllable at the equilibrium point $\boldsymbol{x}^*$. 
\end{Thm}

\begin{Thm} 
Assume that a chemical reaction has a positive equilibrium $\boldsymbol{x}^*$ and that $\rank \boldsymbol{\gamma} = M$. Then the system is controllable in the positive orthant.
\end{Thm}

We prove that chemical reactions in the form of (\ref{eq:kemdiff2}) are controllable almost everywhere without the need to assume the existence of a positive equilibrium; moreover they are controllable in the positive orthant.

\begin{Thm}
A chemical reaction with all reactions rates as control inputs is controllable almost everywhere. Moreover, it is controllable at every point of the positive orthant.
\end{Thm}
\begin{Proof}
The controllability distribution is spanned by the vector fields from the Lie-algebra generated by the vector fields $\boldsymbol{g}_1,\boldsymbol{g}_2,\ldots,\boldsymbol{g}_R$, and the vector fields $[\boldsymbol{f},\boldsymbol{g}_1],[\boldsymbol{f},\boldsymbol{g}_2],\ldots,[\boldsymbol{f},\boldsymbol{g}_R]$, so it is true, that
\begin{equation}
\Span \left\{ \boldsymbol{g}_1,\boldsymbol{g}_2,\ldots,\boldsymbol{g}_R \right\} \subseteq \Delta_C
\end{equation}
and since for every strictly positive (or negative) vector $\boldsymbol{x}$ this subspace is
\begin{equation}
\Span \left\{ \boldsymbol{g}_1,\boldsymbol{g}_2,\ldots,\boldsymbol{g}_R \right\} = \Span \left\{ \boldsymbol{\gamma}(\cdot,1)c_1, \boldsymbol{\gamma}(\cdot,2)c_2,\ldots,\boldsymbol{\gamma}(\cdot,R)c_R\right\}
\end{equation}
with $c_r = \boldsymbol{x}^{\alpha(\cdot,r)}$ being a nonzero constant for a fixed $\boldsymbol{x}>0$ (or $\boldsymbol{x}<0$) for $r=1,2,\ldots,R$, this implies that $\Delta_C$ equals to the image space of $\boldsymbol{\gamma}$, so $\dim \Delta_C = \rank \boldsymbol{\gamma}$. Therefore, if none of the coordinates of $\boldsymbol{x}$ is zero, then the system is controllable.
\end{Proof}
Note that it is possible that $x_m=0$ for some $m \in \{1,2,\ldots,M\}$ and the system is still controllable.

Another result from \cite{farkascontrol} is that if we make a reaction step reversible, and add this reaction step to the original reaction, then under some conditions, the controllability properties of the resulting system does not change. 

Consider a reaction that has a positive equilibrium. Let the input matrix of this reaction linearized in the positive equilibrium be $\boldsymbol{B}$. Take a reaction step, and make it reversible, and add this step to the original reaction. Suppose that the resulting reaction also has a positive equilibrium, and let the input matrix of the linearized model be $\tilde{\boldsymbol{B}}$. Then the rank of these matrices are the same \cite{farkascontrol}, as it is stated in the following theorem.
\begin{Thm}
Let us assume that a given reaction has a positive equilibrium. If we make a reaction step reversible and the obtained reaction also has a positive equilibrium then $\rank \tilde{\boldsymbol{B}} = \rank \boldsymbol{B}$.
\end{Thm} 
We show, that if we make any reaction step reversible, then it does not change the dimension of the controllability distribution regardless of the existence of a positive equilibrium.
\begin{Thm}
If we make any reaction step reversible, then the resulting chemical reaction is controllable almost everywhere; moreover, it is controllable at every point of the positive orthant.
\end{Thm}
\begin{Proof}
Make the $r$th reaction step reversible, i.e. the reaction step vector of the added reaction step is $\boldsymbol{\alpha}(\cdot,r)-\boldsymbol{\beta}(\cdot,r)=-\boldsymbol{\gamma}(\cdot,r)$ and the reactant complex of the new reaction step is defined by the vector $\boldsymbol{\beta}(\cdot,r)$, so the differential equation of the new reaction step is
\begin{equation}
\dot{\boldsymbol{x}}(t) = -\boldsymbol{\gamma}(\cdot,r) \boldsymbol{x}^{\beta(\cdot,r)}(t)u_{R+1}(t) = \boldsymbol{g}_{R+1}(\boldsymbol{x}(t))u_{R+1}(t).
\end{equation}
The $\tilde{\boldsymbol{\gamma}}$ stoichiometric matrix of the extended reaction is
\begin{equation}
\tilde{\boldsymbol{\gamma}} = \left(\begin{array}{ccccccc} \boldsymbol{\gamma}(\cdot,1) & \boldsymbol{\gamma}(\cdot,2) & \ldots & \boldsymbol{\gamma}(\cdot,r) & \ldots & \boldsymbol{\gamma}(\cdot,R) & -\boldsymbol{\gamma}(\cdot,r) \end{array}\right)
\end{equation}
so it immediately follows that $\rank \tilde{\boldsymbol{\gamma}} = \rank \boldsymbol{\gamma}$. The $\tilde{\Delta}_C$ controllability distribution of the extended reaction possesses the property, that 
\begin{equation}
\Span \left\{ \boldsymbol{\gamma}(\cdot,1)c_1,\boldsymbol{\gamma}(\cdot,2)c_2,\ldots,\boldsymbol{\gamma}(\cdot,R)c_R,-\boldsymbol{\gamma}(\cdot,r)c_{R+1}  \right\} \subseteq \tilde{\Delta}_C 
\end{equation}
where the left-hand side is the image space of $\tilde{\boldsymbol{\gamma}}$ if the constants $c_1,c_2,\ldots,c_R$ are not zero.
\end{Proof}

Based on the results of this section, we can conclude that chemical reactions are controllable at almost every point if all the reaction rate coefficients can be varied independently to control the system. Moreover, controllability is not harmed even if we make any reaction step reversible. In the next section we analyze the practically more important case when not all reaction rate coefficients can be used as control inputs.

\section{Controllability of chemical reactions with a reduced number of control parameters}\label{sec:lessinput}

Suppose, that we can not control all the reaction rate coefficients independently, only a subset of them. Suppose without the loss of generality, that there are $K$ number of controllable reaction rate coefficients, and the reaction steps are indexed so that the reaction rate coefficient of the first $K$ reaction steps can be controlled. Let the index set of the controllable reaction rate coefficients be denoted by $I=\{1,2,\ldots,K\}$, while the index set of the remaining reaction rate coefficients be denoted by $N=\{K+1,\ldots,R\}$. In this case $1 \leq K < R$. Let the vector field corresponding to the $r$th reaction step be denoted by $g_r$, i.e. 
\begin{equation}\label{eq:reactionvector}
\boldsymbol{g}_r : x \mapsto \boldsymbol{\gamma}(\cdot,r) \boldsymbol{x} ^{\boldsymbol{\alpha}(\cdot,r)},
\end{equation}
then the differential equation of the reaction with $K$ control inputs is
\begin{equation}\label{eq:kemdiff3}
\dot{\boldsymbol{x}}(t)=\underbrace{\mathop{\sum}_{n \in N} k_n \boldsymbol{g}_n(\boldsymbol{x}(t))}_{\boldsymbol{f}(\boldsymbol{x}(t))} + \mathop{\sum}_{i \in I} \boldsymbol{g}_i(\boldsymbol{x}(t)) u_i(t).
\end{equation}

\begin{Def}
The controllability distribution of the reaction (\ref{eq:kemdiff3}) with $K$ control inputs is the distribution $\Delta_C^{(K)}$, if 
\begin{enumerate}
\item $\forall \boldsymbol{x} \in \mathbb{R}^M$, $\Span \left\{ \boldsymbol{g}_{1}(\boldsymbol{x}),\boldsymbol{g}_{2}(\boldsymbol{x}),\ldots,\boldsymbol{g}_{K}(\boldsymbol{x}) \right\} \subseteq \Delta_C^{(K)}(\boldsymbol{x})$.
\item It is invariant under the vector field $\boldsymbol{f}$, i.e. if $\boldsymbol{\tau} \in \Delta_C^{(K)}$, then $[\boldsymbol{\tau},\boldsymbol{f}] \in \Delta_C^{(K)}$.
\item It is involutive, i.e. if $\boldsymbol{\tau_1},\boldsymbol{\tau_2} \in \Delta_C^{(K)}$, then $[\boldsymbol{\tau_1},\boldsymbol{\tau_2}] \in \Delta_C^{(K)}$.
\end{enumerate}
\end{Def}

The controllability distribution $\Delta_C^{(K)}$ is thus spanned by the vector fields from the Lie-algebra generated by the vector fields $\boldsymbol{g}_1, \boldsymbol{g}_2, \ldots , \boldsymbol{g}_K$ and the vector fields $[\boldsymbol{f},\boldsymbol{g}_1], [\boldsymbol{f},\boldsymbol{g}_2], \ldots , [\boldsymbol{f},\boldsymbol{g}_K] $. Therefore in the following analysis, we will need the Lie-brackets of the vector fields corresponding to reaction steps of the form of (\ref{eq:reactionvector}). We will show, that at each fixed $\boldsymbol{x} \in \mathbb{R}^M$, the Lie-bracket of the vector fields corresponding to two reaction steps is the linear combination of the reaction step vectors. First, we introduce the following notation
\begin{equation}
d_m^r (\boldsymbol{x}) = \left\{ \begin{array}{lr} 0, & \text{ if } \alpha(m,r)=0, \\ \boldsymbol{x}^{\boldsymbol{\alpha}(\cdot,r)-\boldsymbol{e_m}} & \text{otherwise}. \end{array} \right.
\end{equation}
with $\boldsymbol{e_m}$ being the $m$th $M$-dimensional unit vector. Defining the diagonal matrix $D_r(\boldsymbol{x})=\diag \left( d_1^r(\boldsymbol{x}),d_2^r(\boldsymbol{x}),\ldots,d_M^r(\boldsymbol{x}) \right)$, we can write the derivative of the monomial $\boldsymbol{x}^{\boldsymbol{\alpha}(\cdot,r)}$ in the short form
\begin{eqnarray}
(\boldsymbol{x}^{\boldsymbol{\alpha}(\cdot,r)})'&=&\left(\begin{array}{cccc} \alpha(1,r)d_1^r(\boldsymbol{x}) & \alpha(2,r)d_2^r(\boldsymbol{x}) & \ldots & \alpha(M,r)d_M^r(\boldsymbol{x})  \end{array}\right)\nonumber\\
&=&\alpha(\cdot,r)^\top D_r(\boldsymbol{x}).
\end{eqnarray}
Utilizing these notations, we may formalize the Lie-bracket of vector fields corresponding to reaction steps as follows.
\begin{Lem}\label{lem:liebracket1}  
Let $\boldsymbol{g}_i,\boldsymbol{g}_j$ be vector fields in the form of (\ref{eq:reactionvector}). Then their Lie-bracket at the point $\boldsymbol{x} \in \mathbb{R}^M$ is  
\begin{eqnarray}
[\boldsymbol{g}_i,\boldsymbol{g}_j](\boldsymbol{x}) &=& \boldsymbol{\gamma}(\cdot,j)\boldsymbol{\alpha}(\cdot,j)^\top D_j(\boldsymbol{x}) \boldsymbol{\gamma}(\cdot,i)\boldsymbol{x}^{\boldsymbol{\alpha}(\cdot,i)} \\\nonumber
&&-\boldsymbol{\gamma}(\cdot,i)\boldsymbol{\alpha}(\cdot,i)^\top D_i(\boldsymbol{x}) \boldsymbol{\gamma}(\cdot,j)\boldsymbol{x}^{\boldsymbol{\alpha}(\cdot,j)}. 
\end{eqnarray}
\end{Lem}
\begin{Proof}
Since the Lie-bracket of two vector fields is $[\boldsymbol{g}_i,\boldsymbol{g}_j]=\boldsymbol{g}_j'\boldsymbol{g}_i-\boldsymbol{g}_i'\boldsymbol{g}_j$ by (\ref{eq:liebracket}), and the vector fields are of the form of (\ref{eq:reactionvector}), we only have to analyze the derivative of the vector fields. The derivative of the vector field $\boldsymbol{g}_i$ at a point $\boldsymbol{x}\in \mathbb{R}^M$ is
\begin{equation}
\boldsymbol{g}_i'(\boldsymbol{x})=\left( \begin{array}{cccc} \partial_1 g_{i,1}(\boldsymbol{x}) & \partial_2 g_{i,1}(\boldsymbol{x}) & \ldots & \partial_M g_{i,1}(\boldsymbol{x}) \\ \partial_1 g_{i,2}(\boldsymbol{x}) & \partial_2 g_{i,2}(\boldsymbol{x}) & \ldots & \partial_M g_{i,2}(\boldsymbol{x}) \\ \vdots &&&\vdots \\ \partial_1 g_{i,M}(\boldsymbol{x}) & \partial_2 g_{i,M}(\boldsymbol{x}) & \ldots & \partial_M g_{i,M}(\boldsymbol{x}) \end{array}\right)=\left(\begin{array}{c} \left(\gamma(1,i)\boldsymbol{x}^{\boldsymbol{\alpha}(\cdot,i)}\right)'\\ \left(\gamma(2,i)\boldsymbol{x}^{\boldsymbol{\alpha}(\cdot,i)}\right)'\\ \vdots \\ \left(\gamma(M,i)\boldsymbol{x}^{\alpha(\cdot,i)}\right)' \end{array}\right)
\end{equation}
where $g_{i,m}(\boldsymbol{x})=\gamma(m,i)\boldsymbol{x}^{\boldsymbol{\alpha}(\cdot,i)}$. This derivative can be written as
\begin{equation}
\boldsymbol{g}_i'(\boldsymbol{x})=\boldsymbol{\gamma}(\cdot,i)\boldsymbol{\alpha}(\cdot,i)^\top D_i(\boldsymbol{x}),
\end{equation}
where we have utilized the notation introduced before the lemma. Since the vector field at a point $\boldsymbol{x} \in \mathbb{R}^M$ is $\boldsymbol{g}_j(\boldsymbol{x})=\boldsymbol{\gamma}(\cdot,j)\boldsymbol{x}^{\boldsymbol{\alpha}(\cdot,j)}$, the statement of the Lemma follows after substitution into the defining equation (\ref{eq:liebracket}) of the Lie-bracket. 
\end{Proof}

\begin{Lem}\label{lem:liebracketzero}
$\forall \boldsymbol{x} \in \mathbb{R}^M$, $\boldsymbol{\alpha}(\cdot,i)^\top D_i(\boldsymbol{x}) \boldsymbol{\gamma}(\cdot,j)=0$  $\Leftrightarrow$ $\alpha(m,i)\gamma(m,j)=0$, $m = 1,2,\ldots,M$.
\end{Lem}
\begin{Proof}
Expanding the term $D_i$ in the following expression, we get
\begin{eqnarray}
\boldsymbol{\alpha}(\cdot,i)^\top D_i(\boldsymbol{x})\boldsymbol{\gamma}(\cdot,j)&=&\alpha(1,i)\gamma(1,j)d_1^i(\boldsymbol{x}) + \alpha(2,i)\gamma(2,j)d_2^i(\boldsymbol{x}) + \nonumber \\
&& + \ldots + \alpha(M,i)\gamma(M,j)d_M^i(\boldsymbol{x})
\end{eqnarray}
where $d_m^i(\boldsymbol{x})$ is either the zero polynomial, or the derivative of a fixed polynomial with regard to its $m$th variable, for any $m$. However $d_m^i=0$ for all $m \in \left\{1,2,\ldots,M \right\}$ if and only if $\boldsymbol{\alpha}(\cdot,i)=0$, in this case the statement holds trivially. So suppose that $\boldsymbol{\alpha}(\cdot,i)\neq 0$. If $d_m^i=0$ for some $m$, then $\alpha(m,i)=0$ by the definition of $d_m^i$, so $\alpha(m,i)\gamma(m,i)=0$ for that $m$. So we can suppose without the loss of generality that $d_m^i\neq 0$ for any $m \in \left\{1,2,\ldots,M \right\}$. Since the nonzero polynomials $d_m^i$, $m=1,2,\ldots,M$ are all the derivatives of the monomial $\boldsymbol{x}^{\boldsymbol{\alpha}(\cdot,i)}$, it follows that they are linearly independent, i.e.
\begin{equation}
\mathop{\sum}_{m=1}^M \alpha(m,i)\gamma(m,j)d_m^i(\boldsymbol{x})=0 \Leftrightarrow \alpha(m,i)\gamma(m,j)=0,\, m=1,2,\ldots,M
\end{equation}
that concludes the proof.
\end{Proof}

If we use the notation $\odot$ for the elementwise (or Hadamard) product of two vectors, the condition for Lemma \ref{lem:liebracketzero} can be written as $\boldsymbol{\alpha}(\cdot,i) \odot \boldsymbol{\gamma}(\cdot,j)=0$.

\begin{Def}
Let $\kappa_{i,j}: \boldsymbol{x} \mapsto \boldsymbol{\alpha}(\cdot,i)^\top D_i(\boldsymbol{x})\boldsymbol{\gamma}(\cdot,j)$.
\end{Def}
Using this notation, the Lie-bracket of two vector fields of the form (\ref{eq:reactionvector}) is
\begin{equation}\label{eq:liebracketkappa}
[\boldsymbol{g}_i,\boldsymbol{g}_j](\boldsymbol{x})=\kappa_{j,i}(\boldsymbol{x})\boldsymbol{x}^{\boldsymbol{\alpha}(\cdot,i)}\boldsymbol{\gamma}(\cdot,j)-\kappa_{i,j}(\boldsymbol{x}) \boldsymbol{x}^{\boldsymbol{\alpha}(\cdot,j)}\boldsymbol{\gamma}(\cdot,i).
\end{equation}

In the analysis of the Lie-bracket of vector fields corresponding to reaction steps, species that participate in both the reactant and product complex of a reaction step by the same amount do not count. We call this species the direct catalyst of a reaction step.

\begin{Def}
The $m$th species is the direct catalyst of the $r$th reaction step, if the following two (equivalent) statements hold.
\begin{enumerate}
\item $\gamma(m,r)=0$, but $\alpha(m,r)\neq 0$.
\item The $m$th species participates in the reactant and product complex of the $r$th reaction step by the same amount.
\end{enumerate} 
\end{Def}

Note that in the reaction step
\begin{equation}
\ce{X +2Y -> 2X +Y}\nonumber
\end{equation}
neither the species $\ce{X}$ nor the species $\ce{Y}$ is a direct catalyst according to our definition.
\begin{Def}
The $m$th species is the direct catalyst of a reaction, if the following two (equivalent) statements hold.
\begin{enumerate}
\item $\boldsymbol{\gamma}(m,\cdot)=0$.
\item If the $m$th species participates in a reaction step, then it is a direct catalyst of that reaction step.
\end{enumerate}
\end{Def}

\begin{Lem}\label{lem:liebracket2}
The Lie-bracket of two vector fields $\boldsymbol{g}_i,\boldsymbol{g}_j$ corresponding to reaction steps written in the form of (\ref{eq:reactionvector}) at the point $\boldsymbol{x} \in \mathbb{R}^M$ is
\begin{enumerate}
\item $[\boldsymbol{g}_i,\boldsymbol{g}_j](\boldsymbol{x})=0$, if there exists no species in the reaction that participates in both the $i$th and $j$th reaction step or if it participates, it is a direct catalyst of the reaction steps.
\item $[\boldsymbol{g}_i,\boldsymbol{g}_j](\boldsymbol{x})=-\kappa_{i,j}(\boldsymbol{x})\boldsymbol{x}^{\boldsymbol{\alpha}(\cdot,j)}\boldsymbol{\gamma}(\cdot,i)$, if there exists at least one species that participates in the reactant complex of the $i$th reaction step, and participates in the $j$th reaction step, but not as a direct catalyst of the $j$th reaction step, and there exists no species that participates in the reactant complex of the $j$th reaction step but does not participate in the $i$th reaction step except as a direct catalyst.
\item $[\boldsymbol{g}_i,\boldsymbol{g}_j](\boldsymbol{x})=+\kappa_{j,i}(\boldsymbol{x})\boldsymbol{x}^{\boldsymbol{\alpha}(\cdot,i)}\boldsymbol{\gamma}(\cdot,j)$, if there exists at least one species that participates in the reactant complex of the $j$th reaction step, and participates in the $i$th reaction step, but not as a direct catalyst of the $i$th reaction step, and there exists no species that participates in the reactant complex of the $i$th reaction step but does not participate in the $j$th reaction step except as a direct catalyst.
\item $[\boldsymbol{g}_i,\boldsymbol{g}_j](\boldsymbol{x})=\kappa_{j,i}(\boldsymbol{x})\boldsymbol{x}^{\boldsymbol{\alpha}(\cdot,i)}\boldsymbol{\gamma}(\cdot,j)-\kappa_{i,j}(\boldsymbol{x})\boldsymbol{x}^{\boldsymbol{\alpha}(\cdot,j)}\boldsymbol{\gamma}(\cdot,i)$, if there exists a species that participates in the reactant complex of the $i$th reaction step and participates in the $j$th reaction step but not as a direct catalyst, and there exists a (possibly different) species that participates in the reactant complex of the $j$th reaction step and participates in the $i$th reaction step, but not as a direct catalyst. 
\end{enumerate}
\end{Lem}
\begin{Proof}
(1): The condition is equivalent to $\boldsymbol{\alpha}(\cdot,i) \odot \boldsymbol{\gamma}(\cdot,j)=0$ and $\boldsymbol{\alpha}(\cdot,j) \odot \boldsymbol{\gamma}(\cdot,i)=0$, which implies that $[\boldsymbol{g}_i,\boldsymbol{g}_j]=0$ because of Lemma \ref{lem:liebracketzero}.

(2): The statement is equivalent to $\boldsymbol{\alpha}(\cdot,i) \odot \boldsymbol{\gamma}(\cdot,j) \neq 0$ but $\boldsymbol{\alpha}(\cdot,j) \odot \boldsymbol{\gamma}(\cdot,i) = 0$, so by Lemma \ref{lem:liebracketzero} $\kappa_{i,j}\neq 0$ but $\kappa_{j,i} = 0$, that implies $[\boldsymbol{g}_i,\boldsymbol{g}_j](\boldsymbol{x})=-\kappa_{i,j}(\boldsymbol{x})\boldsymbol{x}^{\boldsymbol{\alpha}(\cdot,j)}\boldsymbol{\gamma}(\cdot,i)$ because of (\ref{eq:liebracketkappa}).

(3): The consequence of (2) and the anticommutativity of the Lie-bracket.

(4): The consequence of (2) and (3) and the bilinearity of the Lie-bracket.
\end{Proof}

%

In order to prove our next theorem, we need the following lemma:
\begin{Lem}\label{lem:independentpoly}
Let $P(\boldsymbol{x}) = \left(\begin{array}{cccc} \boldsymbol{f}_1(\boldsymbol{x}) & \boldsymbol{f}_2(\boldsymbol{x}) & \ldots & \boldsymbol{f}_n(\boldsymbol{x})  \end{array}\right)$ be an $n \times n$ matrix whose elements are polynomials of $\boldsymbol{x} \in \mathbb{R}^n$ of finite order. If there exists a point $\boldsymbol{x}^* \in \mathbb{R}^n$ such that $\det P(\boldsymbol{x}^*) \neq 0$, then $\det P \neq 0$ holds almost everywhere in $\mathbb{R}^n$.
\end{Lem}
\begin{Proof}
Since the elements of the matrix are finite order polynomials, the determinant of $P$ is also a finite order polynomial. Since there is a point at which $\det P \neq 0$, this polynomial is not the zero polynomial, thus the set of the roots of the determinant has zero measure in $\mathbb{R}^n$, so $\det P \neq 0$ almost everywhere in $\mathbb{R}^n$. 
\end{Proof}

\begin{Thm}\label{thm:localtoglobal}
If a chemical reaction is controllable at a point $\boldsymbol{x}^* \in \mathbb{R}^M$ with the inputs $K$, then it is controllable almost everywhere.
\end{Thm}
\begin{Proof}
If the reaction is controllable at $\boldsymbol{x}^*$, then $\dim \Delta_C^{(K)} = \rank \boldsymbol{\gamma}$, which means that there are $\rank \boldsymbol{\gamma}$ number of independent vectors spanning $\Delta_C^{(K)}$. The controllability distribution $\Delta_C^{(K)}$ is spanned by the vector fields from the Lie-algebra generated by vector fields of the form of (\ref{eq:reactionvector}), which are polynomials with each monomial being multiplied with one of the reaction step vectors because of Lemmas \ref{lem:liebracket1}--\ref{lem:liebracket2}, thus by collecting the reaction step vectors we get that
\begin{equation}
\Span \left\{ \rho_1(\boldsymbol{x}) \boldsymbol{\gamma}(\cdot,1),\rho_2(\boldsymbol{x}) \boldsymbol{\gamma}(\cdot,2),\ldots,\rho_R(\boldsymbol{x}) \boldsymbol{\gamma}(\cdot,R)  \right\} \subseteq \Delta_C^{(K)}(\boldsymbol{x})
\end{equation}
with $\rho_r,\quad r=1,2,\ldots,R$ being some polynomials. If the reaction is controllable at a point $\boldsymbol{x}^*$, then
\begin{equation}\label{eq:deltackpoly}
\mathop{\sum}_{r=1}^R c_r \rho_r(\boldsymbol{x}^*) \boldsymbol{\gamma}(\cdot,r)=0
\end{equation}
if and only if there are $\rank \boldsymbol{\gamma}$ number of $c_r$s that are zero. Let $\rank \boldsymbol{\gamma} = S$. Suppose without the loss of generality, that the first $S$ component of the vectors $\rho_1(\boldsymbol{x}^*)\boldsymbol{\gamma}(\cdot,1),\rho_2(\boldsymbol{x}^*)\boldsymbol{\gamma}(\cdot,2),\ldots,\rho_S(\boldsymbol{x}^*)\boldsymbol{\gamma}(\cdot,S)$ are linearly independent. Denote the choice of the first $S$ components of the polynomial vector $\rho_i \boldsymbol{\gamma}(\cdot,i)$ by $\rho_i \boldsymbol{\gamma}(\cdot,i)^{[S]}$. Then linear independence can be expressed so that the system of equations defined by
\begin{equation}\label{eq:submat}
\left( \begin{array}{cccc} \rho_1 \boldsymbol{\gamma}(\cdot,1)^{[S]}(\boldsymbol{x}) & \rho_2 \boldsymbol{\gamma}(\cdot,2)^{[S]}(\boldsymbol{x}) & \ldots & \rho_S \boldsymbol{\gamma}(\cdot,S)^{[S]}(\boldsymbol{x}) \end{array}\right) \left(\begin{array}{c} c_1 \\ c_2 \\ \vdots \\ c_S \end{array}\right)=0
\end{equation}
at the point $\boldsymbol{x}^*$ has only the trivial solution. This yields that the determinant of the matrix on the left-hand side of (\ref{eq:submat}) is not zero, and by Lemma \ref{lem:independentpoly}, the determinant of that matrix is not zero almost everywhere in $\mathbb{R}^S$, so the vectors are linearly independent almost everywhere.
\end{Proof}

In the following proof we give an upper estimate for the required number of control inputs.
\begin{Thm}
A reaction is controllable almost everywhere with $\rank \boldsymbol{\gamma}$ number of reaction rate coefficients used as control inputs, i.e. $K = \rank \boldsymbol{\gamma}$.
\end{Thm}
\begin{Proof}
Choose as control inputs the reaction rate coefficients of reaction steps whose reaction step vectors are linearly independent. Since the controllability distribution contains these vector fields corresponding to $\rank \boldsymbol{\gamma}$ number of independent reaction steps, the dimension of the distribution is $\rank \boldsymbol{\gamma}$ almost everywhere.
\end{Proof}

In practice, usually $\rank \boldsymbol{\gamma} < R$, so controlling only $\rank \boldsymbol{\gamma}$ reaction rate coefficients may be much easier than controlling the reaction rate coefficients of every reaction. The question is, whether we can further decrease the number of control inputs.
\begin{Def}
Let the reaction rate coefficients of the first $K$ reaction steps be control inputs. A reaction step is said to be the critical reaction step of the control system with the $K$ inputs, if removing the rate coefficient of that reaction step from the control inputs reduces the dimension of the controllability distribution almost everywhere. 
\end{Def}
Simply saying, if we remove the reaction rate coefficient of a critical reaction step from a controllable system then the resulting system will not be controllable. Identification of critical reaction steps is crucial in finding the required control inputs. 
\begin{Def}
A reaction step is said to be an initializer, if the species of its reactant complex do not participate in any other reaction steps, except if such a species is a direct catalyst in the given step, and the species participating in the reactant complex are not all direct catalysts.
\end{Def}
For example, in the reaction
\begin{align}
\cee{X + U &->[k_1] X + V}\\
\cee{X + Y &->[k_2] X + Z}
\end{align}
both reaction steps are initializers, while in the reaction
\begin{align}
\cee{X + U &->[k_1] U + V}\\
\cee{X + V &->[k_2] X + W}
\end{align}
only the first reaction step is an initializer.

Initializers have key role in the control of chemical reactions, because they are always critical reaction steps.
\begin{Thm}\label{thm:initializer}
An initializer is a critical reaction step for all $K$ input sets that contain the initializer. 
\end{Thm}
\begin{Proof}
Let the $i$th reaction step be an initializer. Due to the definition of the initializer, there exist an $m \in \{1,2,\ldots,M\}$ such that $\gamma(m,i) \neq 0$, but $\gamma(m,j)=0$ for all $j \in \{1,2,\ldots,R\}$, $j\neq i$, so the reaction step vector corresponding to an initializer is linearly independent from all the other reaction step vectors of the chemical reaction.

Suppose indirectly, that there exists a number $K$ such that $\Delta_C^{(K)}$ has dimension  $\rank \gamma$ almost everywhere, but $i \notin I=\{1,2,\ldots,K\}$. We will show, that this can not happen.

Since $i \notin I$, and $\boldsymbol{\gamma}(\cdot,i)$ is linearly independent from other column vectors of $\boldsymbol{\gamma}$, $\boldsymbol{\gamma}(\cdot,i)$ is linearly independent from all $\boldsymbol{g}_k, k\in I$ almost everywhere. Let $\mathcal{G}$ denote the Lie-algebra generated by the vector fields $\boldsymbol{g}_k, k \in I$. Due to Lemma \ref{lem:liebracket2}, this Lie-algebra consists of vector fields that can be written as 
\begin{equation}
\boldsymbol{g} = \mathop{\sum}_{k \in I} \boldsymbol{\gamma}(\cdot,k) \rho_k
\end{equation}
with $\rho_k$ being some polynomials. At almost every $\boldsymbol{x} \in \mathbb{R}^M$, the vector $\boldsymbol{g}(\boldsymbol{x})$ is linearly independent of $\boldsymbol{\gamma}(\cdot,i)$, since it is the linear combination of vectors $\boldsymbol{\gamma}(\cdot,k), k\in I$ that are linearly independent from $\boldsymbol{\gamma}(\cdot,i)$.

Let us examine the vector fields of the form $[\boldsymbol{f},\mathcal{G}]$. Because the Lie-algebra is bilinear, this can be written as
\begin{equation}
[\boldsymbol{f},\mathcal{G}] = \left[k_i \boldsymbol{g}_i,\mathcal{G}\right] + \left[\mathop{\sum}_{n \in N, n \neq i} k_n \boldsymbol{g}_n,\mathcal{G}\right].
\end{equation}
Vector fields of the form $[\mathop{\sum}_{n \in N, n \neq i} k_n \boldsymbol{g}_n,\mathcal{G}]$ are linearly independent of $\boldsymbol{\gamma}(\cdot,i)$ almost everywhere due to the same argument as before, since the vectors $\boldsymbol{g}_n, n \in N, n \neq i$ are linearly independent from $\boldsymbol{\gamma}(\cdot,i)$ almost everywhere.

Albeit $\boldsymbol{g}_i(\boldsymbol{x}) = \boldsymbol{\gamma}(\cdot,i) \boldsymbol{x}^{\boldsymbol{\alpha}(\cdot,i)}$ is linearly dependent of $\boldsymbol{\gamma}(\cdot,i)$, the vector fields $[k_i \boldsymbol{g}_i,\mathcal{G}]$ are also linearly independent of $\boldsymbol{\gamma}(\cdot,i)$, since the $i$th reaction step is an initializer, i.e. there exists no species that is in the reactant complex of the $i$th reaction step, but participates in the other reaction steps except as a direct catalyst, so by Lemma \ref{lem:liebracket2}, the Lie-bracket of these vector fields does not contain a component linearly dependent of $\boldsymbol{\gamma}(\cdot,i)$. 

The Lie-algebra generated by $\mathcal{G}$ and $[\boldsymbol{f},\mathcal{G}]$ is also linearly independent from $\boldsymbol{\gamma}(\cdot,i)$ almost everywhere, since the vector fields in the generators are linearly independent from $\boldsymbol{\gamma}(\cdot,i)$ almost everywhere.

The controllability distribution is formed by vector fields from the Lie-algebra generated by $\mathcal{G}$ and $[\boldsymbol{f},\mathcal{G}]$, but $\boldsymbol{\gamma}(\cdot,i)$ is linearly independent from these vector fields, thus $\dim \Delta_C^{(K)} < \rank \boldsymbol{\gamma}$ almost everywhere, which implies the system is not controllable, that is a contradiction.
\end{Proof}

This theorem states that the reaction rate coefficient of an initializer needs to be controlled in order to control the chemical reaction, so a lower estimate of the required number of control inputs is the number of initializers in a reaction. Note that the number of initializers can be easily found utilizing the $\boldsymbol{\alpha}$ and $\boldsymbol{\gamma}$ matrices. Usually the required number of control inputs is higher than the number of initializers, but there are some special reactions that can be controlled by controlling the reaction rate coefficients of its initializers.

The following are also true about initializers:
\begin{enumerate}
	\item A step of a reversible pair can not be an initializer.
	\item If the reaction is weakly reversible, then there is no initializer, but it is possible to find initializer classes, which we define later.
	\item If any of the species of the reactant complex are zero initially, then the initial step can not start, and the reaction is not controllable at the initial point.
\end{enumerate}

\begin{Def}
A reaction is called a consecutive reaction if it has the following properties.
\begin{enumerate}
\item It has exactly one initializer.
\item If we remove the initializer, the resulting reaction has exactly one initializer (or it is empty), so it has the previous property.
\item If we remove the initializer, the resulting reaction has the previous two properties.
\end{enumerate}
\end{Def}
In consecutive reactions the reaction steps follow each other serially, so they may also be called serial reactions. The examples in Section \ref{sec:examples} will show that our definition of consecutive reactions is more general than the one usually used in chemical kinetics. A consecutive reaction has only one initializer, and the reaction is controllable by controlling the reaction rate coefficient of the initializer only.
\begin{Thm}\label{thm:consecutive}
A consecutive reaction is controllable with only one input, if and only if the input is the reaction rate coefficient of its initializer.
\end{Thm}
\begin{Proof}
If the reaction is controllable with one input, then the input must be the reaction rate coefficient of the initializer due to Theorem \ref{thm:initializer}.

Now we prove, that if the control input is the reaction rate coefficient of the initializer, then the reaction is controllable. Let the initializer of the consecutive reaction have the index $i_1$. If we remove the initializer, we get a consecutive reaction due to the definition of consecutive reactions, or an empty reaction. Let the initializer of the new reaction be the reaction step that has the index $i_2$ in the original reaction. The reactant complex of the $i_2$th reaction step has at least one species that participates in the $i_1$th reaction step and is not a direct catalyst in either reactions, otherwise the original reaction would have the $i_2$th reaction step as an initializer as well. The reaction step with index $i_1$ does not have common (non direct catalyst) species with reaction steps other than the one with index $i_2$, since if it had a common non direct catalyst species with an reaction step with index $j \neq i_2$, it would mean that removing $i_1$ would result in a reaction with both $i_2$ and $j$ being initializers (or none of them), that contradicts the definition of a consecutive reaction. So the reaction step with index $i_1$ has common species only with reaction step with index $i_2$. Remove the reaction steps with indices $i_1$ and $i_2$. The resulting reaction is either empty, or it has one initializer, let this be the reaction step with index $i_3$ in the original reaction. Using similar argument as before, the reaction step with index $i_2$ has a non direct catalyst species that participates in the reactant complex of the reaction step with $i_3$, and the reaction step with index $i_2$ does has common species with reaction steps with indices $i_1$ and $i_3$ only (except for species that are direct catalysts of the $i_2$th reaction step). We find the $i_n$th reaction step by removing the first $n-1$ initializers recursively. Then the $i_n$th reaction step has common species (other than direct catalysts) only with reaction steps with indices $i_{n-1}$ and $i_{n+1}$, and the reactant complex of the reaction step with index $i_n$ has a common non direct catalyst species with the reaction step with index $i_{n-1}$. 

The indices of the reaction steps of the consecutive reaction can be arranged into the set $\{i_1,i_2,\ldots,i_R\}$, and the drift and the control vector fields of the consecutive reaction can be written as
\begin{eqnarray}
\boldsymbol{f} &:& \boldsymbol{x} \mapsto \mathop{\sum}_{r=2}^R \boldsymbol{g}_{i_r} k_r \nonumber\\
\boldsymbol{g} &:& \boldsymbol{x} \mapsto \boldsymbol{g}_{i_1}(\boldsymbol{x}).
\end{eqnarray}
The controllability distribution $\Delta_C^{i_1}$ consists of the vector fields from the Lie-algebra generated by the vector fields $\boldsymbol{g}$ and $[\boldsymbol{f},\boldsymbol{g}]$. The vector field $\boldsymbol{g}$ spans $\boldsymbol{\gamma}(\cdot,i_1)$ almost everywhere because of its definition. The Lie-bracket $[\boldsymbol{f},\boldsymbol{g}]$ can be written as
\begin{equation}
[\boldsymbol{f},\boldsymbol{g}] = \left[ \boldsymbol{g}_{i_2},\boldsymbol{g}_{i_1} \right] + \left[ \mathop{\sum}_{l \neq 2} \boldsymbol{g}_{i_l},\boldsymbol{g}_{i_1} \right]
\end{equation}
where the second term is zero due to Lemma \ref{lem:liebracketzero}, since the reaction step with index $i_1$ has common non direct catalyst species only with the reaction step with index $i_2$. The reaction step with index $i_2$ has a non direct catalyst species in its reactant complex that participates in the $i_1$th reaction step not as a direct catalyst, but there is no species in the reactant complex of the $i_1$th reaction step that participates in the $i_2$th reaction step as a non direct catalyst, so the application of Lemma \ref{lem:liebracket2} results in
 the vector field $[\boldsymbol{f},\boldsymbol{g}]$ at the point $\boldsymbol{x} \in \mathbb{R}^M$ written as 
\begin{equation}
[\boldsymbol{f},\boldsymbol{g}](\boldsymbol{x})=\left[ \boldsymbol{g}_{i_2},\boldsymbol{g}_{i_1} \right](\boldsymbol{x})= \kappa_{i_2,i_1}\boldsymbol{x}^{\boldsymbol{\alpha}(\cdot,i_1)} \boldsymbol{\gamma}(\cdot,i_2),
\end{equation}
so $\boldsymbol{g}$ and $[\boldsymbol{f},\boldsymbol{g}]$ span $\{\boldsymbol{\gamma}(\cdot,i_1),\boldsymbol{\gamma}(\cdot,i_2)\}$ almost everywhere. Due to the definition of $\kappa_{i_2,i_1}$, the vector field $[\boldsymbol{f},\boldsymbol{g}]$ can be interpreted as a vector field corresponding to linear  combination of fictive reaction steps with the same reaction step vector as $\boldsymbol{\gamma}(\cdot,i_2)$, and the reactant complex vectors being element-wise less than or equal to $\boldsymbol{\alpha}(\cdot,i_2)+\boldsymbol{\alpha}(\cdot,i_1)$, so the species in the reactant complex vector of the fictive reaction steps are from the species of the reactant complexes of the reaction steps with indices $i_1$ and $i_2$, that does not participate in any other reaction steps. So the fictive reaction steps related to $[\boldsymbol{f},\boldsymbol{g}]$ has common species with the reactant complex of the $i_3$th reaction step, and the reactant complex of the fictive reaction step related to $[\boldsymbol{f},\boldsymbol{g}]$ has common species with the $i_1$th and $i_2$th reaction step (except for direct catalysts). This implies, that the vector field $[\boldsymbol{f},[\boldsymbol{f},\boldsymbol{g}]]$ can be written as
\begin{equation}
[\boldsymbol{f},[\boldsymbol{f},\boldsymbol{g}]]= \left[\boldsymbol{g}_{i_2},\kappa_{i_2,i_1}\boldsymbol{\gamma}(\cdot,i_2)\right] + \left[\boldsymbol{g}_{i_3},\kappa_{i_2,i_1}\boldsymbol{\gamma}(\cdot,i_2)\right] 
\end{equation}
that can be written at a point $\boldsymbol{x} \in \mathbb{R}^M$ as
\begin{equation}
[\boldsymbol{f},[\boldsymbol{f},\boldsymbol{g}]] (\boldsymbol{x}) = \rho_{1,2}(\boldsymbol{x}) \boldsymbol{\gamma}(\cdot,i_2) + \rho_{2,2}(\boldsymbol{x}) \boldsymbol{\gamma}(\cdot,i_3)
\end{equation}
where $\rho_{1,2}$ and $\rho_{2,2}$ are polynomials in which the monomials divide either $\boldsymbol{x}^{\boldsymbol{\alpha}(\cdot,i_1)}$ or $\boldsymbol{x}^{\boldsymbol{\alpha}(\cdot,i_2)}$. The vector field $[\boldsymbol{f},[\boldsymbol{f},\boldsymbol{g}]]$ can be interpreted as a vector field corresponding to the linear combination of fictive reaction steps with reaction step vectors $\boldsymbol{\gamma}(\cdot,i_2)$ and $\boldsymbol{\gamma}(\cdot,i_3)$ and the species in the reactant complex of the fictive reaction steps are the species from the reactant complexes of the reaction steps with indices $i_1,i_2,i_3$. So the fictive reactions related to the vector field $[\boldsymbol{f},[\boldsymbol{f},\boldsymbol{g}]$ does not have a species in their reactant complex that participates in the reaction steps with index $i_r, r>3$ as a non direct catalyst, but there is a species in the reactant complex of the $i_4$th reaction step that participates in the fictive reactions related to $[\boldsymbol{f},[\boldsymbol{f},\boldsymbol{g}]]$. 
 
Continuing the same argument, we get that application of the Lie-bracket with the vector field $\boldsymbol{f}$ to the vector field $\boldsymbol{g}$ $n-1$ times results in the vector field with the value 
\begin{eqnarray}
\underbrace{[\boldsymbol{f},...[\boldsymbol{f},[\boldsymbol{f}}_{n-1},\boldsymbol{g}]]...](\boldsymbol{x})&=&\rho_{1,n-1}(\boldsymbol{x})\boldsymbol{\gamma}(\cdot,i_2)+\rho_{2,n-1}(\boldsymbol{x})\boldsymbol{\gamma}(\cdot,i_3)\nonumber\\
&&+ \ldots +\rho_{n-1,n-1}(\boldsymbol{x})\boldsymbol{\gamma}(\cdot,i_n)
\end{eqnarray} 
in the point $\boldsymbol{x} \in \mathbb{R}^M$. This shows, that application of the Lie-bracket with vector field $\boldsymbol{f}$ to the vector field $\boldsymbol{g}$ $n-1$ times results in a vector field that spans $\{\boldsymbol{\gamma}({\cdot,i_2}),\boldsymbol{\gamma}(\cdot,i_3),\ldots,\boldsymbol{\gamma}(\cdot,i_n)\}$ almost everywhere if $n \leq R$. This implies that the Lie-algebra generated by the vector fields $\boldsymbol{g}$ and $[\boldsymbol{f},\boldsymbol{g}]$ (and thus the controllability distribution) spans the subspace $\{\boldsymbol{\gamma}(\cdot,i_1),\boldsymbol{\gamma}(\cdot,i_2),\ldots,\boldsymbol{\gamma}(\cdot,i_R)\}$ almost everywhere, so it spans the image space of $\boldsymbol{\gamma}$ almost everywhere, thus the dimension of the controllability distribution is $\rank \boldsymbol{\gamma}$ almost everywhere.
\end{Proof}

It is possible, that a chemical reaction does not have an initializer. However, in this case there is a group of reaction steps such that if they were united into one reaction step, the result would be an initializer. We will call these groups the initializer classes.

\begin{Def}
We say that a group of reaction steps form an initializer class, if the species from the reactant complex of the reaction steps in the initializer class does not participate in the reaction steps that are not in the initializer class, except as direct catalysts.
\end{Def}

\begin{Thm}
There is a critical reaction step in each initializer class.
\end{Thm}
\begin{Proof}
The proof is indirect. Suppose, that the reaction rate coefficients corresponding to the reaction steps in the initializer class are not among the control inputs, and that the chemical reaction is controllable almost everywhere. Consider the fictive reaction that is the linear combination of the reaction steps in the initializer class. This fictive reaction can be used in the analysis, since the drift vector fields of the chemical system in question consists of the linear combination of the vector fields corresponding to the reaction steps on the initializer class and possibly other reaction steps and the Lie-bracket is bilinear. The linear combination of the reaction steps of the initializer class is an initializer by definition, that is a critical reaction step due to Theorem \ref{thm:initializer}, so the chemical reaction is not controllable, that is a contradiction.
\end{Proof}

\section{Examples}\label{sec:examples}

Consider the following consecutive reaction with three species ($M=3$) and two reaction steps ($R=2$)
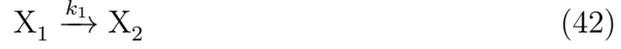
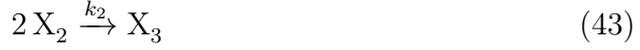
\begin{eqnarray}
\ce{X1 ->[k_1] X2}\\
\ce{2X2 ->[k_2] X3}
\end{eqnarray}
This reaction has one initializer, it is the first reaction step with reaction rate coefficient $k_1$. The matrices $\boldsymbol{\alpha}, \boldsymbol{\beta}, \boldsymbol{\gamma}$ formed by the reactant complex vectors, product complex vectors and reaction step vectors respectively are
\begin{equation}
\begin{array}{ccc} \boldsymbol{\alpha}=\left(\begin{array}{cc} 1 & 0 \\ 0 & 2 \\ 0 & 0 \end{array}\right), & \boldsymbol{\beta}=\left(\begin{array}{cc} 0 & 0 \\ 1 & 0 \\ 0 & 1 \end{array}\right),   & \boldsymbol{\gamma}=\left(\begin{array}{rr} -1 & 0 \\ 1 & -2 \\ 0 & 1 \end{array}\right)\end{array}.
\end{equation}
This is a consecutive reaction, so by Theorem \ref{thm:consecutive} it is controllable with a single input if the control input is the reaction rate coefficient of the initializer, i.e. it is the reaction rate coefficient $k_1$. The control and drift vector fields are
\begin{equation}
\begin{array}{cc} \boldsymbol{g}: \boldsymbol{x}\mapsto \left(\begin{array}{r} -1 \\ 1 \\ 0 \end{array}\right)x_1, & \boldsymbol{f}: \boldsymbol{x}\mapsto k_2\left(\begin{array}{r} 0 \\ -2 \\ 1 \end{array}\right)x_2^2  \end{array},
\end{equation} 
and the differential equation of the chemical reaction is
\begin{eqnarray}
\dot{\boldsymbol{x}}(t) &=& \boldsymbol{f}(\boldsymbol{x}(t)) + \boldsymbol{g}(\boldsymbol{x}(t))u(t) \nonumber\\
&=& k_2\left(\begin{array}{r} 0 \\ -2 \\ 1 \end{array}\right)x_2^2(t)+\left(\begin{array}{r} -1 \\ 1 \\ 0 \end{array}\right)x_1(t) u(t).
\end{eqnarray}
where $\boldsymbol{x}=(x_1,x_2,x_3)^\top$ and $x_1,x_2$ and $x_3$ is the concentration of the species $\ce{X1},\ce{X2}$ and $\ce{X3}$, respectively. The controllability distribution is spanned by the vector fields from the Lie-algebra generated by $\boldsymbol{g}$ and $[\boldsymbol{f},\boldsymbol{g}]$. The vector field $[\boldsymbol{f},\boldsymbol{g}]$ at the point $\boldsymbol{x} \in \mathbb{R}^3$ is
\begin{eqnarray}
[\boldsymbol{f},\boldsymbol{g}](\boldsymbol{x})&=&(\boldsymbol{g}'\boldsymbol{f})(\boldsymbol{x})-(\boldsymbol{f}'\boldsymbol{g})(\boldsymbol{x}) \nonumber\\
&=& \left(\begin{array}{rrr}-1 & 0 & 0 \\ 1 &0 & 0 \\ 0 & 0 & 0\end{array}\right)\left(\begin{array}{r}0 \\-2 \\ 1\end{array}\right)k_2 x_2^2 - k_2\left(\begin{array}{rrr}0 & 0 & 0 \\ 0 & -4x_2 & 0 \\ 0 & 2x_2 & 0\end{array}\right)\left(\begin{array}{r}-1\\1\\0\end{array}\right)x_1\nonumber\\
&=&-2k_2x_1x_2\left(\begin{array}{r} 0 \\ -2 \\ 1 \end{array}\right)=-2k_2x_1x_2\boldsymbol{\gamma}(\cdot,2).
\end{eqnarray}
The vector fields $\boldsymbol{g}$ and $[\boldsymbol{f},\boldsymbol{g}]$ are linearly independent almost everywhere, so the dimension of the controllability distribution is at least two almost everywhere. Since the dimension of the stoichiometric space is $\rank \boldsymbol{\gamma}=2$, the reaction is controllable almost everywhere.

Now we determine the set of null measure on which the system is not controllable. The system is trivially not controllable on the line $x_1=0$, since in this case $\boldsymbol{g}(\boldsymbol{x})=0$. Now suppose, that $x_1 \neq 0$. The vector field $[\boldsymbol{f},\boldsymbol{g}]$ is zero at $x_2=0$, however there may be other vector fields in the Lie-algebra generated by the vector fields $\boldsymbol{g}$ and $[\boldsymbol{f},\boldsymbol{g}]$ that are not zero at $x_2=0$. 

The vector field $[\boldsymbol{g},[\boldsymbol{f},\boldsymbol{g}]]$ at a point $\boldsymbol{x} \in \mathbb{R}^3$ is
\begin{equation}
[\boldsymbol{g},[\boldsymbol{f},\boldsymbol{g}]](\boldsymbol{x}) = ([\boldsymbol{f},\boldsymbol{g}]'\boldsymbol{g}-\boldsymbol{g}'[\boldsymbol{f},\boldsymbol{g}])(\boldsymbol{x})=2k_2x_1(x_2-x_1)\left(\begin{array}{r} 0\\-2\\1\end{array}\right),
\end{equation} 
that is not zero, if $x_2=0$, however it is zero, if $x_1=x_2$. So if $x_1 \neq 0$ and $x_2 \neq 0$, then the stoichiometric space is spanned by the vector fields $\boldsymbol{g}$ and $[\boldsymbol{f},\boldsymbol{g}]$, and if $x_1 \neq 0$ and $x_2=0$, then the stoichiometric space is spanned by the vector fields $\boldsymbol{g}$ and $[\boldsymbol{g},[\boldsymbol{f},\boldsymbol{g}]]$. This implies that the reaction is controllable everywhere except on the line $x_1=0$.

Consider the following consecutive reaction that consists of three reaction steps with four species:
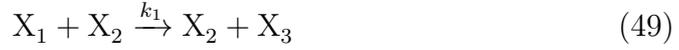
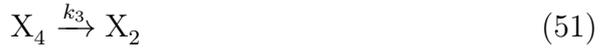
\begin{align}
\cee{X1 + X2 &->[k_1] X2 + X3}\\
\cee{X3 &->[k_2] X4 }\\
\cee{X4 &->[k_3] X2 }
\end{align} 
This reaction has one initializer, that is the reaction step with the reaction rate coefficient $k_1$. The matrices $\boldsymbol{\alpha}, \boldsymbol{\beta}, \boldsymbol{\gamma}$ formed by the reactant complex vectors, product complex vectors and reaction step vectors respectively are
\begin{equation}
\begin{array}{ccc} \boldsymbol{\alpha}=\left(\begin{array}{ccc} 1 & 0 & 0 \\ 1 & 0 & 0 \\ 0 & 1 & 0 \\ 0 & 0 & 1 \end{array}\right), & \boldsymbol{\beta}=\left(\begin{array}{ccc} 0 & 0 & 0 \\ 1 & 0 & 1 \\ 1 & 0 & 0 \\ 0 & 1 & 0  \end{array}\right),   & \boldsymbol{\gamma}=\left(\begin{array}{rrr} -1 & 0 & 0 \\ 0 & 0 & 1 \\ 1 & -1 & 0 \\ 0 & 1 & -1 \end{array}\right)\end{array}.
\end{equation}
Let the control input be the $k_1$ reaction rate coefficient of the first reaction, and let the concentrations of the species $\ce{X1}$, $\ce{X2}$, $\ce{X3}$ and $\ce{X4}$ be denoted by $x_1$, $x_2$, $x_3$ and $x_4$ respectively. Then the control and drift vector fields are
\begin{equation}
\begin{array}{cc} \boldsymbol{g}: \boldsymbol{x}\mapsto \left(\begin{array}{r} -1 \\ 0\\ 1 \\ 0 \end{array}\right)x_1x_2, & \boldsymbol{f}: \boldsymbol{x}\mapsto k_2\left(\begin{array}{r} 0 \\ 0 \\ -1 \\ 1 \end{array}\right)x_3+k_3\left(\begin{array}{r} 0 \\ 1 \\ 0 \\ -1 \end{array}\right)x_4  \end{array},
\end{equation} 
and the differential equation of the chemical reaction system is
\begin{eqnarray}
\dot{\boldsymbol{x}}(t) &=& \boldsymbol{f}(\boldsymbol{x}(t)) + \boldsymbol{g}(\boldsymbol{x}(t))u(t) \nonumber\\
&=& \left(\begin{array}{c} 0 \\ k_3x_4(t) \\ -k_2x_3(t) \\ k_2x_3(t)-k_3x_4(t) \end{array}\right)+\left(\begin{array}{r} -1 \\ 0\\ 1 \\ 0 \end{array}\right)x_1(t)x_2(t) u(t).
\end{eqnarray}

We calculate the vector fields $[\boldsymbol{f},\boldsymbol{g}]$ and $[\boldsymbol{f},[\boldsymbol{f},\boldsymbol{g}]]$. The vector field $[\boldsymbol{f},\boldsymbol{g}]$ at the point $\boldsymbol{x} \in \mathbb{R}^4$ is
\begin{eqnarray}
[\boldsymbol{f},\boldsymbol{g}](\boldsymbol{x})&=&(\boldsymbol{g}'\boldsymbol{f})(\boldsymbol{x})-(\boldsymbol{f}'\boldsymbol{g})(\boldsymbol{x}) \nonumber\\
&=& \left(\begin{array}{rrrr}-x_2 & -x_1 & 0 & 0 \\ 0  & 0 & 0 & 0 \\ x_2 & x_1 & 0 & 0 \\ 0 & 0 &0 &0 \end{array}\right)\left(\begin{array}{c} 0\\ k_3x_4 \\ -k_2x_3 \\ k_2x_3-k_3x_4\end{array}\right)\nonumber\\
&& - \left(\begin{array}{rrrr}0 & 0 & 0 & 0 \\ 0&  0 & 0 & k_3 \\ 0 & 0 & -k_2 & 0 \\ 0 & 0 & k_2 & -k_3 \end{array}\right)\left(\begin{array}{r}-1\\0\\1\\0\end{array}\right)x_1x_2\nonumber\\
&=&x_1 \left( \begin{array}{c}-k_3x_4 \\ 0 \\ k_2x_2+k_3x_4\\-k_2x_2 \end{array}\right).
\end{eqnarray}
The vector field $[\boldsymbol{f},[\boldsymbol{f},\boldsymbol{g}]]$ at the point $\boldsymbol{x} \in \mathbb{R}^4$ is
\begin{eqnarray}
[\boldsymbol{f},[\boldsymbol{f},\boldsymbol{g}]](\boldsymbol{x})&=&([\boldsymbol{f},\boldsymbol{g}]'\boldsymbol{f}-\boldsymbol{f}'[\boldsymbol{f},\boldsymbol{g}])(\boldsymbol{x})\nonumber\\
&=&x_1\left(\begin{array}{c}-k_2k_3x_3+k_3^2x_4\\k_2k_3x_2\\k_2^2x_2-k_3^2x_4+k_2k_3x_3+2k_2k_3x_4\\-k_2(k_2x_2+k_3x_2+2k_3x_4) \end{array}\right).
\end{eqnarray}
The vector fields $\boldsymbol{g}$, $[\boldsymbol{f},\boldsymbol{g}]$ and $[\boldsymbol{f},[\boldsymbol{f},\boldsymbol{g}]]$ are linearly independent almost everywhere, and they are in the Lie-algebra generated by the vector fields $\boldsymbol{g}$ and $[\boldsymbol{f},\boldsymbol{g}]$, so the controllability distribution contains these vector fields, and thus the dimension of the controllability distribution is three almost everywhere. Since $\rank \boldsymbol{\gamma} = 3$, the chemical reaction system is controllable almost everywhere.

Now we determine the set of points in which the system is not controllable. Since $\boldsymbol{g}(\boldsymbol{x})=0$ if $x_1=0$ or $x_2=0$, it follows that the system is not controllable in the plane $x_1=x_2=0$. Suppose now, that $x_1 \neq 0$ and $x_2 \neq 0$. It can be easily verified that the vector fields $\boldsymbol{g},[\boldsymbol{f},\boldsymbol{g}]$ and $[\boldsymbol{f},[\boldsymbol{f},\boldsymbol{g}]]$ are linearly independent in this case. So the system is controllable in every $\boldsymbol{x} \in \mathbb{R}^4$ except if $x_1=0$ or $x_2=0$.

\section{Discussion}\label{sec:discuss}

We have shown that if all the reaction rate coefficients of a chemical reaction are control inputs, then the chemical reaction is controllable almost everywhere, and it is controllable everywhere in the positive orthant. We have shown, that making a reaction step reversible and adding it to the original reaction does not affect controllability.

In chemical engineering practice, reducing the number of control inputs is of great interest. We have defined critical reaction steps of input sets, that must remain in the input set to keep controllability. However, the state of being a critical reaction step depends on the current input set. We have defined initializer reaction steps, and showed that they are always critical reaction steps independent of the input set. Thus the minimal number of control inputs can not be less than the number of initializer reaction steps. We have defined consecutive reactions, and proved that they can be controlled with a single input.

It is possible, that a reaction does not have an initializer reaction step, however it is still possible to find an initializer class of the reaction steps. We have proved, that there is always at least one critical reaction step in every initializer class. So the minimal number of required control inputs is not less than the sum of the number of initializer reaction steps and initializer classes.

Most of the results are specific for chemical reactions, and the theorems and definitions are expressed using chemical concepts like e.g. the species participating in the reactant and product complexes. However, Theorem \ref{thm:localtoglobal} only uses the fact that chemical reactions are polynomial systems, so this theorem applies for any system whose dynamics is described by a differential equation with its right-hand side being a polynomial of the state variables, and being linear in the control inputs. So in the case of input affine polynomial systems, local controllability implies global controllability. 

It is a nontrivial question how the inputs are to be chosen in order to achieve a given final state, we shall return to this --- important from the practical point of view --- question later.

In the future, we plan to investigate the observability of chemical reactions using the same approach. By using the Lie-algebra of the smooth vector fields, we hope to find global results for observability to generalize the local results that can be found e.g. in \cite{farkasobserv,ZsHorvath}.

\subsection*{Acknowledgements}

The present work has partially been supported by the 
Hungarian National Scientific Foundation, No. 84060, and the CMST COST Action CM1404 Chemistry of Smart Energy Carriers and Technologies (SMARTCATS).

\bibliography{drexlertoth}
\bibliographystyle{abbrv}
\end{document}